\documentclass[preprint]{elsarticle}

\usepackage{graphicx}
\usepackage{lineno}
\usepackage{amssymb}
\usepackage{mathrsfs}
\usepackage{mhchem}
\usepackage{subcaption}
\usepackage{epstopdf}
\usepackage{float}
\usepackage{bbold}
\usepackage{amsmath}
\usepackage{color} 
\usepackage{ulem} 
\usepackage{algorithmicx}
\usepackage{algorithm}
\usepackage{algpseudocode}

\usepackage{caption}
\usepackage{hyperref}
\usepackage{amsfonts}
\usepackage{graphicx}

\newcommand{\R}{\mathbb{R}}

\renewcommand{\u}{\mathbf{u}}

\renewcommand{\c}{\mathbf{c}}

\numberwithin{theorem}{section}
\usepackage{amsopn}

\date{\today}

\begin{document}

\begin{frontmatter}
\title{Adaptive Time-stepping Schemes for the Solution of the Poisson-Nernst-Planck Equations}

\author{David Yan}
\address{Department of Electrical and Computer Engineering, University of Toronto}

\author{M.C. Pugh\corref{mycorrespondingauthor}}
\address{Department of Mathematics, University of Toronto, \\
40 St George St, Toronto, ON  M5S 2E4, Canada}
\cortext[mycorrespondingauthor]{Corresponding author}
\ead{mpugh@math.utoronto.ca}

\author{F.P. Dawson}
\address{Department of Electrical and Computer Engineering, University of Toronto}


\begin{abstract}
The Poisson-Nernst-Planck equations with generalized Frumkin-Butler-Volmer boundary conditions (PNP-FBV) describe ion transport with Faradaic reactions, and have applications in a number of fields. In this article, we develop an adaptive time-stepping scheme for the solution of the PNP-FBV system based on two time-stepping methods: a fully-implicit (VSBDF2) method, and a semi-implicit (VSSBDF2) method. We present simulations under both current and voltage boundary conditions and demonstrate the ability to simulate a large range of parameters, including any value of the singular perturbation parameter $\epsilon$.  When
the underlying dynamics is one that would have the solutions converge to a steady-state solution, we observe that the adaptive time-stepper based on the VSSBDF2 method produces solutions that ``nearly'' converge to the steady-state solution and that, simultaneously, the time-step sizes stabilize to a limiting size $dt_\infty$. In the companion to this article \cite{YPD_Part2}, we linearize the SBDF2 scheme about the steady-state solution, and demonstrate that the linearized scheme is conditionally stable, and that this conditional stability is the cause of the adaptive time-stepper's behaviour.  
While the adaptive time-stepper based on the fully-implicit (VSBDF2) method is not subject to such time-step stability restrictions, the required nonlinear solve incurs additional computational cost. We profile both methods to identify regimes of the perturbation parameter $\epsilon$ where one method is favourable over the other.
The matlab code used in this work can be found at \url{https://github.com/daveboat/vssimex_pnp} .
\end{abstract}



\end{frontmatter}


\section{Introduction}\label{introduction}

The Poisson-Nernst-Planck (PNP) equations describe the transport of charged species subject to diffusion and 
electromigration.   They have wide applicability in electrochemistry, and have been used to model a  
number of different systems, including porous 
media \cite{Biesheuvel2010PRE, Biesheuvel2011, Biesheuvel2012, Peters2016},
microelectrodes \cite{Streeter2008, Compton2011},
ion-exchange membranes \cite{Dydek2013, Nikonenko2010},
electrokinetic phenomena \cite{Yaroshchuk2012, Bazant2010, Bazant2009},
ionic liquids \cite{Bazant2011, Kornyshev2007},
electrochemical thin films \cite{Bazant2005, Chu2005, Biesheuvel2009},
fuel cells \cite{BiesheuvelFranco2009},
supercapacitors \cite{Lee2014},
and more. 

The one-dimensional, nondimensionalized  PNP equations 
for a media with $2$ mobile species is
\begin{align}
\label{concentration_nondim}
\frac{\partial c_\pm}{\partial t} &= -\frac{\partial}{\partial x}\left[-\frac{\partial c_\pm}{\partial x} - z_\pm \, c_\pm \, \frac{\partial \phi}{\partial x}\right], \qquad t>0, \, x \in (0,1),\\
\label{poisson_nondim}
-\epsilon^2\frac{\partial^2 \phi}{\partial x^2} &= \frac{1}{2}\left(z_+ \, c_+ + z_- \, c_-\right), \qquad \qquad \qquad x \in (0,1),
\end{align}
where $c_\pm$ and $z_\pm$ are the concentration and charge number of the positive/negative ion, $\phi$ is the potential and $\epsilon$ is the ratio of the Debye screening length to the inter-electrode width $L$.   This width is used
in the nondimensionalization of the original modelling equations \cite{Yan2017}; the
domain $(0,L)$ is rescaled to $(0,1)$.
We consider a model in which the anion and cation have a single charge ($z_\pm = \pm 1$) and
the anion has no charge-transfer reactions at the electrode: $c_-$ has no-flux boundary conditions:
\begin{equation}
\label{no_flux_BCs}
- \left( - \frac{\partial  c_-}{\partial x} +  \, c_- \, \frac{\partial \phi}{\partial x}\right)\bigg |_{x=0} \mkern-12mu = \left(- \frac{\partial  c_-}{\partial x}+  \, c_- \, \frac{\partial \phi}{\partial x} \right)\bigg |_{x=1}  = 0.
\end{equation}
The cation is assumed to have a reaction at the electrodes involving the transfer of one electron; this is modelled
using  
generalized Frumkin-Butler-Volmer (FBV) boundary conditions:
\begin{align}
\label{bv_nondim_L}
- \left( - \frac{\partial  c_+}{\partial x} -  \, c_+ \, \frac{\partial \phi}{\partial x}\right)\bigg |_{x=0} \mkern-12mu &= F(t) := 4k_{c,a} \,
c_+(0,t) \, e^{- 0.5 \; \Delta \phi_\text{left}} - 4 \, j_{ r,a} \, e^{0.5 \; \Delta \phi_\text{left}}, \\
\label{bv_nondim_R}
\left(- \frac{\partial  c_+}{\partial x} -  \, c_+ \, \frac{\partial \phi}{\partial x} \right)\bigg |_{x=1} \mkern-12mu &= G(t) := 4k_{c,c} \,
c_+(1,t) \, e^{-0.5 \; \Delta \phi_\text{right}} - 4 \, j_{r,c} \, e^{0.5 \; \Delta \phi_\text{right}},
\end{align}
where $k_{c,a}$, $k_{c,c}$, $j_{r,a}$, and $j_{r,c}$ are reaction rate parameters; the second letter in the subscripts ($a$ and $c$) refer to the anode and cathode, respectively. Equations \eqref{bv_nondim_L}--\eqref{bv_nondim_R} model the electrodeposition reaction
 \begin{equation}
\label{reaction}
\ce{C^+ + e^- <=> M }
\end{equation}
where M represents the electrode material.

The Stern layer is a compact layer of charge that occurs in the electrolyte next to an electrode surface
\cite{Bazant2013ACR, Soestbergen2012}; $\lambda_S$ denotes the effective
width of this layer.  In equations \eqref{bv_nondim_L}--\eqref{bv_nondim_R}, $\Delta
\phi_\text{left}$ and $\Delta \phi_\text{right}$ refer to the
potential differences across the Stern layers that occur at the anode
and cathode respectively.  Specifically,
\begin{equation} \label{Delta_phi_is}
\Delta \phi_\text{left} = \phi_\text{anode} - \phi(0,t) = -\phi(0,t), 
\; \Delta \phi_\text{right} = \phi_\text{cathode} - \phi(1,t) = v(t)-\phi(1,t)
\end{equation}
where the potential at the anode has been set to zero and $v(t)$
denotes the potential at the cathode.  

The Poisson
equation \eqref{poisson_nondim} uses a mixed (or Robin) boundary condition \cite{Bazant2005,
  Chu2005, Biesheuvel2009},
\begin{align}
\label{phi_bc_nondim_L}
- \epsilon \, \delta \; \frac{\partial \phi}{\partial x}\bigg |_{x=0} & = \Delta \phi_\text{left} := - \phi(0,t), \\
\label{phi_bc_nondim_R}
+ \epsilon \, \delta \; \frac{\partial \phi}{\partial x}\bigg |_{x=1} & =\Delta \phi_\text{right} := v(t) - \phi(1,t),
\end{align}
where $\delta=\lambda_S/L$.  Finally, there is an ODE which ensures conservation of electrical
current at the electrode \cite{Moya1995, Soestbergen2010},
\begin{equation}
\label{current_conservation_nondim}
-\frac{\epsilon^2}{2} \, \frac{d \; }{dt} \phi_x(1,t)
=j_\text{ext}(t) - \left[k_{c,c} \, c_+\left(1,t\right) \, e^{-0.5 \; \Delta \phi_\text{right}} - j_{r,c} \, e^{0.5 \; \Delta \phi_\text{right}}\right],
\end{equation}
where $j_\text{ext}(t)$ is the current through the device.  We refer
to the PNP equations with the generalized Frumkin-Butler-Volmer
boundary conditions as the PNP-FBV system.

Considering more than two charged species or reactions involving the transfer
of more than one electron is a straight-forward generalization \cite{Bard2001,Newman2012}.
If the model is extended to include adsorption effects at the
electrodes, there would be an additional ODE describing the dynamics
of the fraction of surface coverage for each electrode
\cite{YanThesis}. If, in addition, the model included temperature and
heat transport, there would be another PDE for the temperature
\cite{YanThesis}. In this work however, we limit ourselves to just the
PNP-FBV system.

The device is operated in two regimes --- either the current or the
voltage at the cathode is externally controlled.  If the voltage at
the cathode, $v(t)$, is externally controlled then the the PNP-FBV
system \eqref{concentration_nondim}--\eqref{poisson_nondim} with
boundary conditions \eqref{no_flux_BCs}--\eqref{bv_nondim_R} and
\eqref{Delta_phi_is}--\eqref{phi_bc_nondim_R} are numerically solved,
determining $c_\pm$ and $\phi$.  The current is found {\it a
  postiori} using equation \eqref{current_conservation_nondim}.  If
the current, $j_\text{ext}(t)$, is externally controlled, then equation
\eqref{current_conservation_nondim} is part of the PNP-FBV system and
the ODE is numerically solved along with the PDEs,
determining $c_\pm$, $\phi$, and $\phi_x(1,t)$ simultaneously.
The voltage $v(t)$ is then found {\it a postiori}. \\

\noindent {\bf Some prior numerical approaches to the PNP and PNP-FBV systems}

\noindent
Though many workers in the field have approximated solutions to the
PNP-FBV system using asymptotic methods \cite{Newman1965, Smyrl1966}, a
numerical method is needed to obtain a full solution. A key
mathematical aspect of the PNP-FBV system is that the parameter
$\epsilon$ acts as a singular perturbation to the system and results
in the formation of boundary layers \cite{Bazant2004}. This makes
numerical simulation of the PNP-FBV system especially challenging.

An accurate, efficient solver for the PNP-FBV equations must take into account the solution's spatial and temporal characteristics. First, 
the bulk ($10 \epsilon \lesssim x \lesssim 1 - 10 \epsilon$) is approximately {\it electroneutral}   \cite{Bockris_Reddy} and so 
the concentrations and potential are linear or nearly linear there, whereas they are nonlinear, and can change significantly over short distances,
near
the electrodes ($0 < x \lesssim 10 \epsilon$ and $1-10 \epsilon \lesssim x <1$). 
A non-uniform mesh which distributes more mesh points near $x=0, 1$ is 
therefore an economical way to discretize the spatial domain. Secondly, many 
electrochemical systems of interest are subject to sudden changes in forcing, 
separated by periods with constant or no forcing. This, combined with the transient 
dynamics driven by initial conditions, gives the problem more than one time-scale, 
motivating the need for specialized time-stepping methods.


One of the first contributions to the field of time-stepping methods
for the Poisson-Nernst-Planck equations, was by Cohen and Cooley
\cite{Cohen1965}, who used an explicit time-stepping method with
predictor-corrector time-step refinement to solve the electroneutral
equations with constant current. The next contribution was by Sandifer and Buck
\cite{Sandifer1975}, who solved a system of time-independent Nernst-Planck equations with constant current and used an
implicit, iterative method to step the displacement field equations. Brumleve and
Buck \cite{Brumleve1978} solved the full PNP equations with
Chang-Jaffe boundary conditions using Backward Euler time-stepping. In
Brumleve and Buck's method, time-steps were allowed to be variable,
but were not adjusted during each step; the method given in the
original publication still sees use in more modern work
\cite{Sokalski2003}. Murphy et. al. \cite{Murphy1993} solved the full PNP-FBV 
system with
adaptive steps by treating the discretized parabolic-elliptic system
as a differential-algebraic system of ODE's and algebraic equations,
and used a variable-order Gear's method \cite{Gear1971}. It
is also worth mentioning here the work of Scharfetter and Gummel
\cite{Scharfetter1969}, who gave a numerical method to solve the
drift-diffusion equations (an analogue of the PNP equations in
semiconductor physics) using Crank-Nicolson time-stepping.

One recent work where the time-dependent PNP equations (with FBV) are solved
is Soestbergen, Biesheuvel and Bazant \cite{Soestbergen2010}. They used the commercial finite element software COMSOL, which uses BDF and the generalized-$\alpha$
method \cite{Erlicher2002}. Another is
Britz and Strutwolf \cite{Britz2014}, who simulate a liquid junction using BDF with constant time-steps.
Britz \cite{Britz2005} also outlines various explicit
and implicit multistep methods in his book on computational solutions
to the PNP equations. \\

\noindent {\bf Our numerical approaches and results of simulations}

\noindent
Our present work differs from previous work
(notably Murphy et. al.) by using a splitting method for the
parabolic-elliptic system
\eqref{concentration_nondim}--\eqref{poisson_nondim} 
and differs from other modern
computational methods \cite{Compton2014} by controlling time-steps
adaptively. In this article, we develop and test two adaptive time-steppers for the
PNP-FBV system with error control, one of which is semi-implicit (also called implicit-explicit) and one of which is fully-implicit.  The adaptivity allows the 
time-stepper to automatically detect changes in time scales and vary the step size accordingly.  

The semi-implicit adaptive time-stepper is based on a second-order variable step-size, semi-implicit, backward
differentiation formula (Variable Step-size Semi-implicit Backwards Differentiation Formula, or VSSBDF2 \cite{Wang2008}).  Its constant step-size counterpart
(SBDF2) is used as a time-stepping scheme in a variety of fields in
science, engineering and computational mathematics. Axelsson
et. al. \cite{Axelsson2014}, for example, used SBDF2 with constant
steps to solve the time-dependent Navier-Stokes equations. Lecoanet
et. al. \cite{Lecoanet2016} used a similar method to model combustion
equations in stars, and Linde, Persson and Sydow \cite{Linde2009} used
it to solve the Black-Scholes equation in computational finance.
An example of an adaptive SBDF2 time-stepping scheme can be found in Rosam, Jimack and Mullis \cite{Rosam2007}, who use it to study a problem in binary alloy solidification. 

The fully-implicit adaptive time-stepper is based on a second-order Variable Step-size, fully-implicit, Backward
Differentiation Formula (VSBDF2).  Its constant step-size counterpart
(BDF2) has the advantage of being unconditionally stable when
applied to linear systems that have asymptotically stable steady states. 
An example of VSBDF2 is Eckert et. al. \cite{Eckert2004}, who use it to model
electroplasticity; their scheme is fully-implicit and they do not report any kind of instability or time
step stability restriction.
However, for nonlinear problems, this numerical stability comes at the cost of needing to perform a
computationally expensive nonlinear solve at each time-step.



We apply the (fully-implicit) VSBDF2 and (semi-implicit) VSSBDF2 
adaptive time-steppers to the PNP-FBV system with time-dependent imposed forcing (voltage or current).  Both time-steppers behave as desired: the 
time-steps refine by orders of magnitude in response to fast changes 
in the forcing and coarsen if the forcing is changing slowly (or not at all).
See Figure \ref{dt_vs_time}.

\begin{figure}[htb!]
\centering
\includegraphics[width=0.47\linewidth]{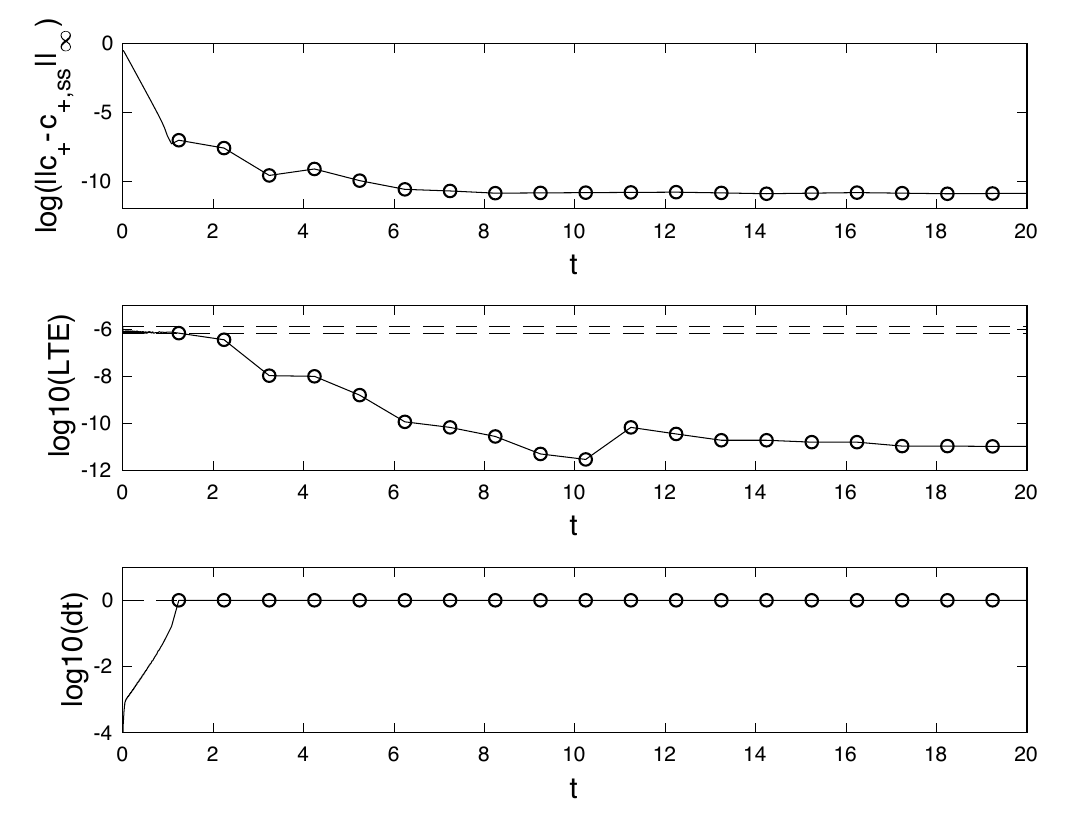}
\includegraphics[width=0.47\linewidth]{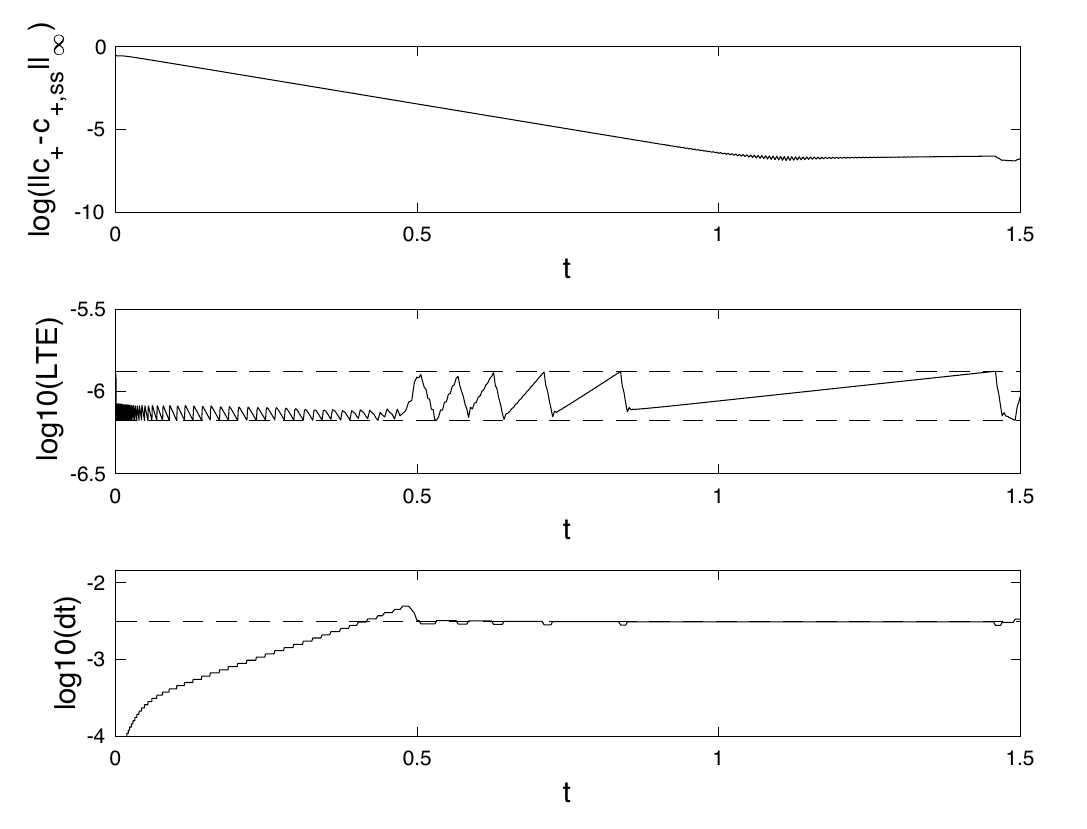}
\caption{PNP-FBV system \eqref{concentration_nondim}--\eqref{phi_bc_nondim_R} 
with constant imposed voltage $v(t) = 2$, $\epsilon = .05$, and all other physical parameters set to $1$.
The initial data is $c_\pm(x,0) = 1 + .1 \, \sin(2 \pi x)$ and $\phi_x(1,0) = 0$.
A uniform mesh is used: $N = 90$ and $L_1 = L_2 = L_3 = 1/3$.
In the three plots on the left, the (fully-implicit) VSBDF2 adaptive time-stepper is used
and, in the three plots on the right, the (semi-implicit) VSSBDF2 adaptive time-stepper is used. In the left plots, open circles are to denote the time steps, after the time-step size $dt$ reaches $dt_\text{max}$. 
The behaviour shown in the right plots between $t=1.25$ and $t=1.5$ is similar to that observed up to $t=20$.
\underline{Top plots:} The logarithm of the maximum deviation of the computed solution $c_+$ from the steady state $c_{+, ss}$ at each moment in 
time: $\log(\|\c_+^n - \c_{+,ss}\|_\infty)$.  Deviations of $\c_-^n$ and $\pmb{\phi}^n$ from the corresponding steady
state profiles behave similarly, as do their time derivatives as approximated
using \eqref{sbdf2}.  \underline{Middle plots:} The logarithm of the approximate local truncation error, \eqref{final_error}.
The dashed lines in both plots indicate $\log(tol \pm range)$. \underline{Bottom plots:} The logarithm of the time-step size, $dt$, plotted versus time. The dashed line in the left plot indicates the maximum step-size $dt_\text{max}$. The dashed line in the right plot indicates the
stability restriction $dt^* = .0031$
computed using the linear stability analysis presented in Section 5
of the companion paper \cite{YPD_Part2}. 
}
\label{Full_Model_tmax_effect}
\end{figure}

However, there are fundamental differences between
the (fully-implicit) VSBDF2 and (semi-implicit) VSSBDF2 
adaptive time-steppers in the constant-forcing regime. In this regime, the
steady-state solution is asymptotically stable.
We find that the (fully-implicit) VSBDF2 adaptive time-stepper takes
larger and larger time steps until reaching 
a user-specified maximum time-step $dt_\text{max}$ and that solutions converge to the steady-state
solution.  See Figure \ref{Full_Model_tmax_effect}.

In contrast, we find that 
the (semi-implicit) VSSBDF2 adaptive time-stepper coarsens its time steps
but they stabilize at a limiting time-step size, $dt_\infty$, which is smaller 
than $dt_\text{max}$.  In addition,
we observe that the solution gets close to, but fails to converge to, the stable steady state. 
%
In the companion article \cite{YPD_Part2} we perform a linear stability 
analysis of the SBDF2 scheme about the steady-state solution, and demonstrate
that the scheme is conditionally stable with a step size stability restriction $dt^*$.  We
demonstrate that it is this underlying stability restriction that causes the
VSSBDF2 adaptive time-stepper to have a limiting time-step size, and present numerical evidence that $dt_\infty = dt^*$.  We find that the stability restriction does not depend significantly 
on the mesh width.  Specifically, $dt_\infty$ does not
go to zero as the mesh width goes to zero.

We find that $dt_\infty \approx \epsilon^2$ for small values of the singular perturbation parameter $\epsilon$; as a result the (semi-implicit) VSSBDF2
adaptive time-stepper will take small step sizes for small values of $\epsilon$.
Although the time steps can be taken larger when using the fully-implicit VSBDF2 adaptive time-stepper, each time step takes longer due to the nonlinear solve required. We compare performance times for both adaptive time-steppers, and identify regimes of the perturbation parameter $\epsilon$ where the runtime for each method is favourable.

Previous work on
time-dependent solutions to the PNP-FBV system \cite{Soestbergen2010, Olesen2010} took comparably large values,
$\epsilon \geq 10^{-3}$, whereas in principle we
are able to simulate with any value of $\epsilon$, allowing exploration of a wider range of operating regimes. For example, the adaptivity of the 
time-stepper allowed us to 
numerically explore certain standard experimental protocols, such as
linear sweep voltammetry, in which a time-dependent voltage or current
is imposed \cite{Yan2017}. 

Separately, we explore how the temporal discretization of the boundary conditions affects
the accuracy of the (semi-implicit) VSSBDF2 scheme.  
We find that 
the ``natural'' semi-implicit approach to the discretizing boundary conditions --- discretizing the boundary conditions by treating the linear term implicitly and extrapolating the nonlinear terms forward in time --- results in the loss of second-order accuracy. In order to retain second-order accuracy, we use a `ghost point' approach to the boundary conditions, where we assume that the PDE holds at the boundaries. 

While the methods we present in the current work are applied to solving the 
PNP-FBV equations, they are more generally relevant to systems of coupled stiff
parabolic-elliptic equations with nonlinear boundary conditions. Also, the stability
analysis in the companion paper is applicable to any multistep semi-implicit scheme.

\subsection{Structure of the article.}
This article is structured as follows. Subsection \ref{time_stepping} presents
the temporal discretization and Subsections 
\ref{spatial_discretization} and \ref{boundary_conditions}
present the spatial discretization.
Subsection \ref{time_stepping_error_control} present the adaptive
stepping and error control algorithms. 
Subsection \ref{splitting_method} discusses how to apply the time-stepping scheme
to the  PNP-FBV system.
In Section \ref{numerical_tests_PNP}, we apply
the scheme to the PNP-FBV system and study its performance including the effect
of the singular perturbation parameter $\epsilon$.
In Appendix \ref{app:extrapolation} we present
the derivation of the local truncation error formula.  

The matlab code used in this work can be found at \url{https://github.com/daveboat/vssimex_pnp} .

\section{The Numerical Method}

\subsection{Time-Stepping}\label{time_stepping}
Multistep schemes such as backwards differentiation formulae, Adams-Bashforth, and Crank-Nicolson methods have long been applied in computational fluid mechanics to time step advection-diffusion systems (see Chapter 4.4 of \cite{Peyret2002}) where both the diffusion and advection terms are linear. 

Semi-implicit, or implicit-explicit schemes, are useful when the system contains both linear stiff terms and nonlinear terms which are difficult to handle using implicit methods. 
Consider the ODE
$u'= f(u)+g(u)$
where $f(u)$ is a nonlinear term and $g(u)$ is a stiff linear term.  
Given $u^{n-1}$ at time
$t^{n-1}=t^n-dt_\text{old}$ and $u^n$ at time $t^n$, 
$u^{n+1}$ at time $t^{n+1}=t^n+dt_\text{now}$ is determined via
\begin{equation}
\text{SBDF2:} \hspace{.2in} \frac{1}{dt}\left(\frac{3}{2}u^{n+1} - 2 \, u^{n} + \frac{1}{2}u^{n-1}\right)
=2 \, f(u^n)- f(u^{n-1}) + g(u^{n+1}), \label{sbdf2}
\end{equation}
where the superscript notation denotes time levels: $u^n$
approximates $u(t^n)$ (see, for example, \cite{Crouzeix_1980,Ascher1995}).

Our VSSBDF2 adaptive time-stepper is based on a second-order variable step-size semi-implicit backwards differencing formula, introduced by Wang and Ruuth \cite{Wang2008}, as a generalization of the  SBDF2 scheme:
\begin{align}\notag
&\text{VSSBDF2:} \hspace{.2in} \frac{1}{dt_\text{now}}\left(\frac{1+2\omega}{1+\omega}u^{n+1} - (1+\omega)u^{n} + \frac{\omega^2}{1+\omega}u^{n-1}\right)
\\
&\hspace{2in}=(1+\omega)f(u^n)-\omega f(u^{n-1}) + g(u^{n+1}), \label{vssbdf2}
\end{align}
where $\omega= dt_\text{now}/dt_\text{old}$.  We use this scheme for
\eqref{concentration_nondim} with term $c_{\pm,xx}$ implicit (like
$g$) and the term $z_\pm \, (c_\pm \, \phi_x)_x$ explicit (like $f$).  We chose
this scheme because the PDEs 
(\ref{concentration_nondim}) are advection diffusion equations; this
scheme was shown to have favorable stability properties compared to
the other IMEX schemes based on numerical experiments on Burger's
equation \cite{Wang2008}.
A one-step semi-implicit scheme is used for the first time-step
\begin{equation} \label{one_step_IMEX}
\frac{1}{dt}\left(u^{1} - u^{0} \right)
= f(u^0) + g(u^{1}).
\end{equation}

To compare the semi-implicit schemes to fully-implicit schemes, we also consider the second-order backwards differencing formula, BDF2,
\begin{equation}
\text{BDF2:} \hspace{.2in} \frac{1}{dt}\left(\frac{3}{3}u^{n+1} - 2 \, u^{n} + \frac{1}{2}u^{n-1}\right)
= g(u^{n+1}). \label{bdf2}
\end{equation}
Our VSBDF2 adaptive time-stepper is based on the second-order variable step-size backwards differencing formula which 
is a generalization of the  BDF2 scheme:
\begin{equation}
\text{VSBDF2:} \hspace{.2in} \frac{1}{dt_\text{now}}\left(\frac{1+2\omega}{1+\omega}u^{n+1} - (1+\omega)u^{n} + \frac{\omega^2}{1+\omega}u^{n-1}\right)
=g(u^{n+1}), \label{vsbdf2}
\end{equation}
where $\omega= dt_\text{now}/dt_\text{old}$.  
Here, $g$ represents $c_{\pm,xx} + z_\pm \, (c_\pm \, \phi_x)_x$.  
Backward Euler is used for the first time-step
\begin{equation} \label{backward_euler}
\frac{1}{dt}\left(u^{1} - u^{0} \right)
= g(u^{1}).
\end{equation}

\subsection{Adaptive Time-Stepping and Error Control}
\label{time_stepping_error_control}
If one is using a non-adaptive time-stepper, the time steps are chosen before computation begins, and are usually constant. On the other hand, an adaptive 
time-stepper chooses the next time-step size $dt_\text{now}$ based on the previously computed solutions and time-step sizes. Specifically, we choose $dt_\text{now}$ by approximating the local truncation error (LTE) $\epsilon_c^{n+1}$, and seeking a $dt_\text{now}$ such that $\epsilon_c^{n+1} \in (tol-range, tol+range)$ \footnote{In practice, we have used a relationship such as $tol=10^{-6}$ and $range=tol/3$.}.

\begin{figure}[htb!]
\centering
\includegraphics[width=0.3\linewidth]{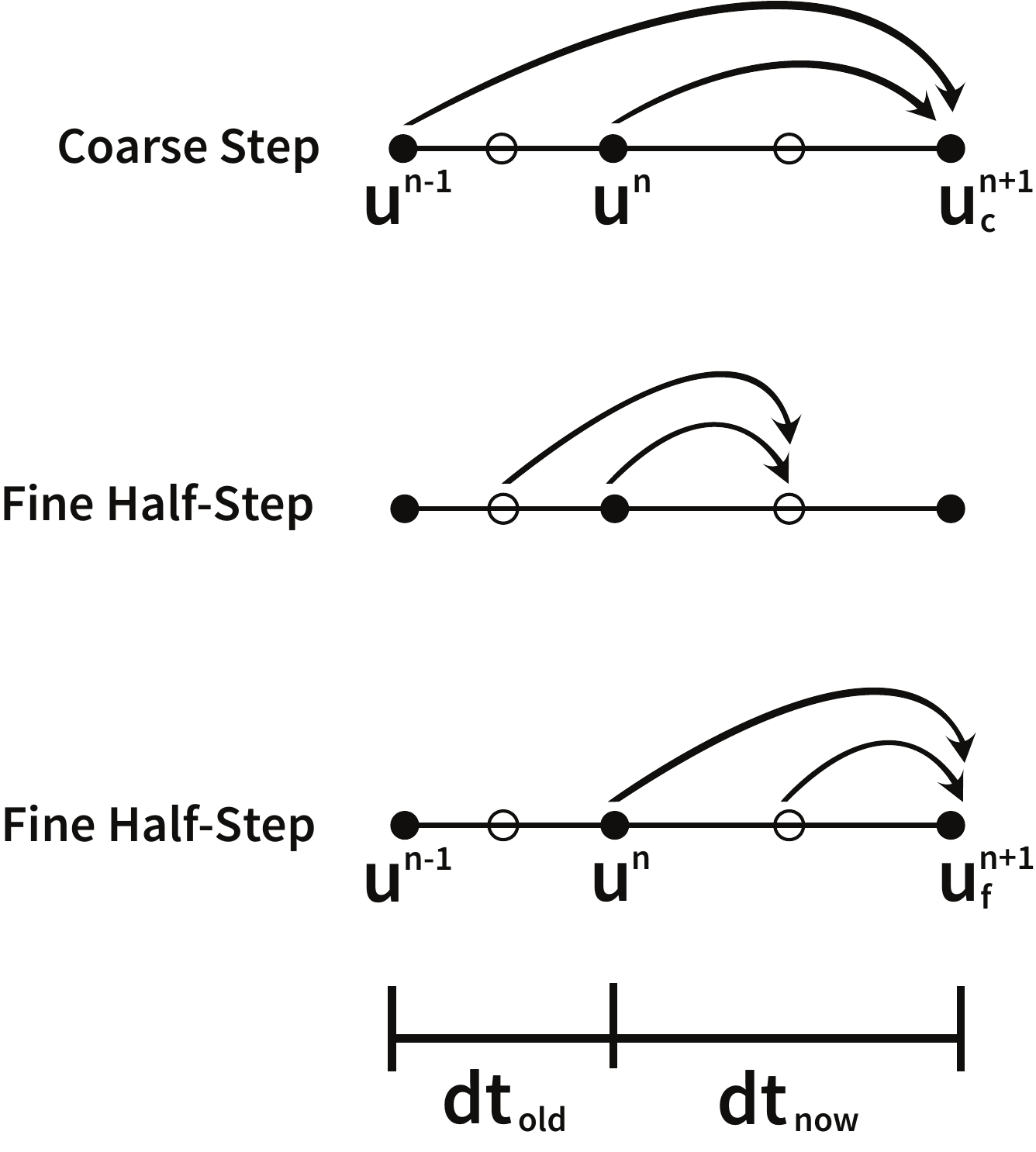}
\caption{Coarse/fine time-stepping scheme. First a coarse time-step uses $u^{n-1}$ and $u^{n}$ to create $u_c^{n+1}$. Then a fine time-step uses $u_f^{n-\frac{1}{2}}$ and $u^{n}$ to create $u_f^{n+\frac{1}{2}}$ and another 
fine time-step uses $u^n$ and $u_f^{n+\frac{1}{2}}$ to create $u^{n+1}_f$
In the diagram, open circles indicate values from fine half-steps.}
\label{timestepping}
\end{figure}

We use a coarse-fine refinement strategy
to approximate the local truncation error: see Figure \ref{timestepping}.
A single ``coarse" time step, resulting in $u^{n+1}_c$, and two ``fine" half time 
steps, resulting in $u^{n+1}_f$, are used to approximate the LTE (see Appendix \ref{app:extrapolation}):
\begin{equation}
\label{final_error}
\epsilon^{n+1}_c = \frac{8\left(dt_\text{old}+dt_\text{now}\right)}{7dt_\text{old}+5dt_\text{now}}\left(u^{n+1}_c-u^{n+1}_f\right)
 \approx u_c^{n+1}-u(t^{n+1}).
\end{equation}
We use this approximation of $\epsilon_c$ when testing whether to use $dt_\text{now}$ or to coarsen or refine it. If $dt_\text{now}$ has been accepted, we can then use $u_f^{n+1}$ and $u_c^{n+1}$ to construct an improved approximation $u^{n+1}$ which has a smaller truncation error. Specifically, $u^{n+1}$ can be
taken as a linear combination of the form
$
u^{n+1}=\alpha u^{n+1}_c + \beta u^{n+1}_f
$
with coefficients
\begin{equation}
\label{final_extrapolation}
\alpha=-\frac{dt_\text{old}+3 \, dt_\text{now}}{7 \, dt_\text{old}+5 \, dt_\text{now}}
\quad \mbox{and} \quad
\beta=8\, \frac{dt_\text{old}+dt_\text{now}}{7 \, dt_\text{old}+5 \, dt_\text{now}}.
\end{equation}
The local truncation error for $u^{n+1}$ is one order higher than the local truncation 
errors for $u^{n+1}_c$ and $u^{n+1}_f$ (see Appendix \ref{app:extrapolation}). Note 
that if $dt_\text{now}=dt_\text{old}$, then \eqref{final_extrapolation} reduces to the 
standard Richardson extrapolation formula for second-order schemes.

As discussed in the companion paper \cite{YPD_Part2}, the use of
Richardson extrapolation affects the stability of the VSSBDF2 time-stepper.
For this reason, Richardson extrapolation was not used in the simulations presented
in this article.

\begin{algorithm}
	\caption{Adaptive time-stepping scheme for a single time step} \label{adaptive_algorithm}
	\begin{algorithmic}
		\State $i \gets 0$ \Comment{Reset loop counter for this time step}
		\State $dt_\text{now} \gets dt_\text{old}$ \Comment{Initial guess at $dt$ for this time step}
		\State $u_c^{n+1} \gets$ TimeStep($dt_\text{now}$) \Comment{Coarse step, TimeStep() using Eq. \eqref{vssbdf2} or \eqref{vsbdf2}}
		\State $u_f^{n+1} \gets$ TimeStep($dt_\text{now}/2$) \Comment{Fine step}
		\State $\epsilon_c^i \gets$ Error($dt_\text{now}$, $dt_\text{old}$, $u_c^{n+1}$, $u_f^{n+1}$) \Comment{Error() from equation \eqref{final_error}}
		\While{abs($\epsilon_c^i$ - tol) $>$ range} \Comment{Loop until the error is acceptable}
			\If{$i \geq i_\text{max}$} \Comment{Enforce maximum iterations}
			\State $u_c^{n+1} \gets$ TimeStep($dt_\text{min}$)
			\State $u_f^{n+1} \gets$ TimeStep($dt_\text{min}/2$)
			\State \textbf{break}
			\EndIf
			\State $dt_\text{now} \gets \text{min}\left(\text{max}\left(\left(\frac{\text{tol}}{\epsilon_c^i}\right)^{1/p}, \eta_\text{min}\right), \, \eta_\text{max}\right)dt_\text{now}$ \Comment{Update $dt$}
			\If{$dt_\text{now} > dt_\text{max}$} \Comment{Enforce maximum time step}
				\State $u_c^{n+1} \gets$ TimeStep($dt_\text{max}$)
				\State $u_f^{n+1} \gets$ TimeStep($dt_\text{max}/2$)
				\State \textbf{break}
			\EndIf
			\State $u_c^{n+1} \gets$ TimeStep($dt_\text{now}$)
			\State $u_f^{n+1} \gets$ TimeStep($dt_\text{now}/2$)
			\State $\epsilon_c^{i+1} \gets$ Error($dt_\text{now}$, $dt_\text{old}$, $u_c^{n+1}$, $u_f^{n+1}$) \Comment{Update error estimate}
			\State $i \gets i+1$
		\EndWhile
		\If{Richardson extrapolation is used}
		\State $u^{n+1} \gets \alpha u_c^{n+1} + \beta u_f^{n+1}$ \Comment{$\alpha$ and $\beta$ are defined in equation \eqref{final_extrapolation}}
		\Else
		\State $u^{n+1} \gets u_c^{n+1}$
		\EndIf
	\end{algorithmic}
\end{algorithm}
Algorithm \ref{adaptive_algorithm} shows our complete time-stepping scheme. The time step update formula $$dt_\text{now} \gets \min\left(\max\left(\left(\frac{tol}{\epsilon^i_c}\right)^{1/p},\eta_\text{min}\right),\eta_\text{max}\right) dt_\text{now}$$
in Algorithm \ref{adaptive_algorithm} is a natural extension of refinement strategies for single step schemes, where $p$ is the order of the LTE. It is based on the idea that the LTE satisfies $\epsilon^i\approx C\left(dt^i\right)^p$, so setting $dt^{i+1}=\left(\frac{tol}{\epsilon^i}\right)^{1/p}dt^i$ will bring $\epsilon^{i+1}$ closer to $tol$. In order to promote convergence, $\eta_\text{min}$ and $\eta_\text{max}$ place limits on how much $dt_\text{now}$ is allowed to change each iteration.
This type of error control strategy is discussed in Chapter II.4 of Hairer, Norsett and Wanner \cite{Hairer2009}.   For the VSBDF2 and VSSBDF2 schemes, $p=3$.
Unless noted otherwise, for the simulations presented here we used
$tol = 10^{-6}$, $range = tol/3$, $\eta_\text{max} = 1.1$, $\eta_\text{min} = .9$,
$dt_\text{max} = 1$ and $dt_\text{min} = 10^{-8}$.

Finally, since we are using a two-step time-stepping scheme, for the first time-step, we use either a one-step IMEX scheme (equation \eqref{one_step_IMEX}) or backward Euler (equation \eqref{backward_euler}), both of which have LTE $\sim O(dt^2)$, along with the error 
estimate $\epsilon^1_c = (4/3) (u^{1}_c-u^{1}_f )$, the time-step update formula
with $p=2$, and the extrapolation formula $u^{1}=2u^1_f-u^1_c$.

Above, we described the time-stepping for a single ODE according to the steps presented above. Given the system of ODEs that would arise from spatially discretizing the conservation laws \eqref{concentration_nondim}, we use the same approach but define the local truncation error using the norm of $\mathbf u$. In our simulations, we chose the $l^2$ norm. 

\subsection{Spatial Discretization}\label{spatial_discretization}

The geometry is divided into a non-uniform mesh, with $x \in [0,1]$ and parameterized via
\begin{equation}
\label{param1}
x(s): [0,1]\rightarrow [0,1],
\end{equation}
so that $x_i:=x(s_i)$, $0=x_1<x_2<...<x_N=1$ and $dx(i):=x(s_{i+1})-x(s_i)$ where $s_i=(i-1)\, ds$, $i=1,2,...,N$ and $ds=1/(N-1)$. The function $x(s)$ may, for example, be piecewise linear 
with smaller slopes near the endpoints $x=0,1$ and a large slope around $x=1/2$.  This would result in a piecewise uniform mesh that is finer near the endpoints. Alternatively, one might use a logistic function, for example, to create a mesh with smoothly varying nonuniformity. In practice, we used  a piecewise uniform mesh with three or five regions of uniform mesh.

To present the spatial discretization of the parabolic PDE, we refer to a generic continuity equation $u_t=-\left(J(u,u_x,x)\right)_x$. For reference, if equation \eqref{concentration_nondim} were written in this general form, we would have $u = c_\pm$ and $J(u,u_x,x) = - c_{\pm,x} - z_\pm \, c_\pm \, \phi_x$.
The function $u(\cdot,t)$ on the interval is approximated by a vector $\mathbf {u}(t) \in \mathbb R^N$ with $u_i(t) \approx u(x_i,t)$.
At internal nodes, the flux is approximated using a second-order center-differencing scheme
\begin{equation}
\label{model_disc2}
\frac{du_i}{dt}\approx (u_t)\big |_{x_i}=\left(J(u,u_x,x)\right)\big |_{x_i} \approx \frac{J(u,u_x,x)\big |_{x_{i+1/2}} - J(u,u_x,x)\big |_{x_{i-1/2}}}{x_{i+1/2}-x_{i-1/2}},
\end{equation}
where $x_{i+1/2}$ is the midpoint of $[x_i,x_{i+1}]$ and the approximations
\begin{equation}
\label{model_disc3}
u_{i\pm 1/2} \approx \frac{u_{i\pm 1}+u_i}{2}, \quad 
u_x\big |_{x_{i - 1/2}} \approx \frac{u_i - u_{i - 1} }{ dx_{i-1} }, 
\quad  \text{and} \quad 
u_x\big |_{x_{i+ 1/2}} \approx \frac{ u_{i+ 1} - u_{i}}{ dx_{i}}
\end{equation}
are used. At the boundary nodes, we use a three-node, second-order approximation for $u_x$. For example, the approximation at the left hand boundary is
\begin{equation}
\label{lhs_derivative}
u_x\big |_{x=0}\approx -\frac{2dx_1+dx_2}{dx_1\left(dx_1+dx_2\right)}u_1 + \frac{dx_1+dx_2}{dx_1dx_2}u_2 - \frac{dx_1}{dx_2\left(dx_1+dx_2\right)}u_3.
\end{equation}

For the spatial discretization of the elliptic PDE (\ref{poisson_nondim}) in the interior, we again use a three point center-differencing scheme
\begin{equation}
\frac{2u_{i+1}}{dx_i \, (dx_i+dx_{i-1})} - \frac{2u_i}{dx_i \, dx_{i-1}}
+\frac{2u_{i-1}}{dx_{i-1} \, (dx_i+dx_{i-1})} \approx u_{xx}\big |_{x_i} = f(x_i)
\end{equation}
where $u = - \phi/\epsilon^2$ and $f = (z_+ \, c_+ + z_- \, c_-)/2$.
On the boundaries $x=0,1$, we use a left or right-handed three point stencil
to approximate the first derivatives in the boundary conditions
(\ref{phi_bc_nondim_L})--(\ref{phi_bc_nondim_R}).

Bazant and coworkers (i.e. in \cite{Chu2005} and \cite{Biesheuvel2009}), used a Chebyshev pseudospectral spatial discretization in their work, where we use a finite difference discretization. We also tested a Chebyshev spectral version of the code using the \texttt{chebfun} package \cite{Trefethen2013, Platte2008} and found that, while the spectral and finite difference codes gave nearly identical results, the finite difference code ran orders of magnitude faster than the spectral code when using the same time-stepping scheme.

\subsection{Boundary Conditions}
\label{boundary_conditions}

For the semi-implicit schemes \eqref{sbdf2} and \eqref{vssbdf2}, a natural first approach to the boundary conditions on the flux of $c_+$,
\eqref{bv_nondim_L}--\eqref{bv_nondim_R}, would be to 
handle the linear term of the flux implicitly and to 
extrapolate both the flux constraint function and the nonlinear term in
the flux forward to time $t^{n+1}$.  We refer to this approach as the
``direct" method of handling the boundary conditions.  For the boundary condition \eqref{bv_nondim_R} 
the ``direct'' method yields 
\begin{equation}\label{bc2_direct}
-(c_{+,x})_N^{n+1} - (1+\omega)( c_+ \phi_x )_N^{n}+\omega(c_+ \phi_x)_N^{n-1} = (1+\omega) \, G^{n}-\omega \, G^{n-1}
\end{equation}
where 
\begin{equation}
G^n = 4k_{c,c} \,
(c_+)_N^n \, e^{-0.5 \; \Delta \phi_\text{right}} - 4 \, j_{r,c} \, e^{0.5 \; \Delta \phi_\text{right}}
\quad \mbox{with} \quad
\Delta \phi_\text{right} =  v(t^n)-\phi_N^n
\end{equation}
and $G^{n-1}$ is defined analogously.
The boundary conditions \eqref{bv_nondim_L} and \eqref{no_flux_BCs} are discretized analogously.
The equations for the discretized boundary conditions and the time-stepping equations at the interior
nodes are then simultaneously solved for $\c_\pm^{n+1}$. 

We find that using the ``direct'' method for boundary conditions
yields the following undesirable properties.  First, the simulation is not second-order accurate in time when the
SBDF2 time-stepper \eqref{sbdf2}
is used: see Table \ref{convergencetest}. Second, we find that we can not take the tolerance to be
arbitrarily small in the VSSBDF2 adaptive time-stepper \cite{YanThesis,YPD_arXiv}.  For these reasons, we use a different approach for the boundary
conditions, which we refer to as the ``ghost point" method (see Section
1.4 of Thomas \cite{Thomas1995}), since it assumes the PDEs
\eqref{concentration_nondim} hold at the end points and uses the same
time-stepping scheme as the one used at the internal nodes.
Applying the PDEs at $x_1$ and $x_N$ requires the flux at neighbouring points; in a true ghost
point method this flux would be located at ``ghost points'' outside
the computational domain: $x_0$ and $x_{N+1}$.  Instead, we use the flux {\it at} $x_1$ and
$x_N$; there are no points outside the computational domain. For example, the spatially discretized PDE for $c_+$ at the left boundary point, $x=0$, is the ODE
\begin{equation}
\label{ghostpoints_left}
\frac{\partial c_+}{\partial t} \bigg |_{x=0} = \frac{(c_{+,x}+c_+ \, \phi_x)\big |_{x= dx_1/2} -F(t)}{dx_1/2}.
\end{equation}
For the semi-implicit time-stepping schemes,  \eqref{sbdf2} and \eqref{vssbdf2}, we
treat the $c_{+,x}$ term implicitly (like $g$) and
extrapolate the $c_+ \, \phi_x$ and $F(t)$ terms (like $f$).

For the fully-implicit schemes \eqref{bdf2}--\eqref{vsbdf2}, we used the ghost point method as well.  The terms on the
right-hand side of \eqref{ghostpoints_left} were all handled implicitly (like $g$).

We find that using the ghost point method for the boundary conditions
yields a method that is second-order accurate in time, when the (constant-time-step)
SBDF2 time-stepper \eqref{sbdf2} or BDF2 time-stepper \eqref{bdf2} are used; see Table \ref{convergencetest}.
Further, we find that the
tolerance can be taken to be as small as we desired, near to round-off, 
in the VSSBDF2 and VSBDF2 adaptive time-steppers \cite{YanThesis,YPD_arXiv}.

\subsection{Time-stepping the Parabolic-Elliptic System}
\label{splitting_method}

The PNP-FBV system has 
2 parabolic PDEs \eqref{concentration_nondim}, one elliptic
PDE \eqref{poisson_nondim}, and may have one ODE 
\eqref{current_conservation_nondim}.  We now describe how to
time-step this system in a way that the global truncation error is
$O(dt^2)$ when we take constant time-steps
($dt_\text{old}=dt_\text{now}=dt$).  We present the case when a
current is imposed and so the ODE \eqref{current_conservation_nondim}
needs to be time-stepped along with the PDEs.

\subsubsection{Using the Semi-Implicit Schemes SBDF2 and VSSBDF2}
We use a splitting scheme.  
The potential $\phi$ appears in the term $(c_\pm \, \phi_x)_x$ in \eqref{concentration_nondim}
and in the boundary conditions \eqref{no_flux_BCs}--\eqref{bv_nondim_R}.
When the semi-implicit schemes SBDF2 \eqref{sbdf2} or VSSBDF2 \eqref{vssbdf2} are used
to time-step the parabolic PDEs, these terms involving $\phi$ are extrapolated forward.  This allows
the new concentrations to be 
computed using the old potentials.  The new potential can be then be computed using the new concentrations.  
\begin{figure}[htb!]
\centering
\includegraphics[width=0.4\linewidth]{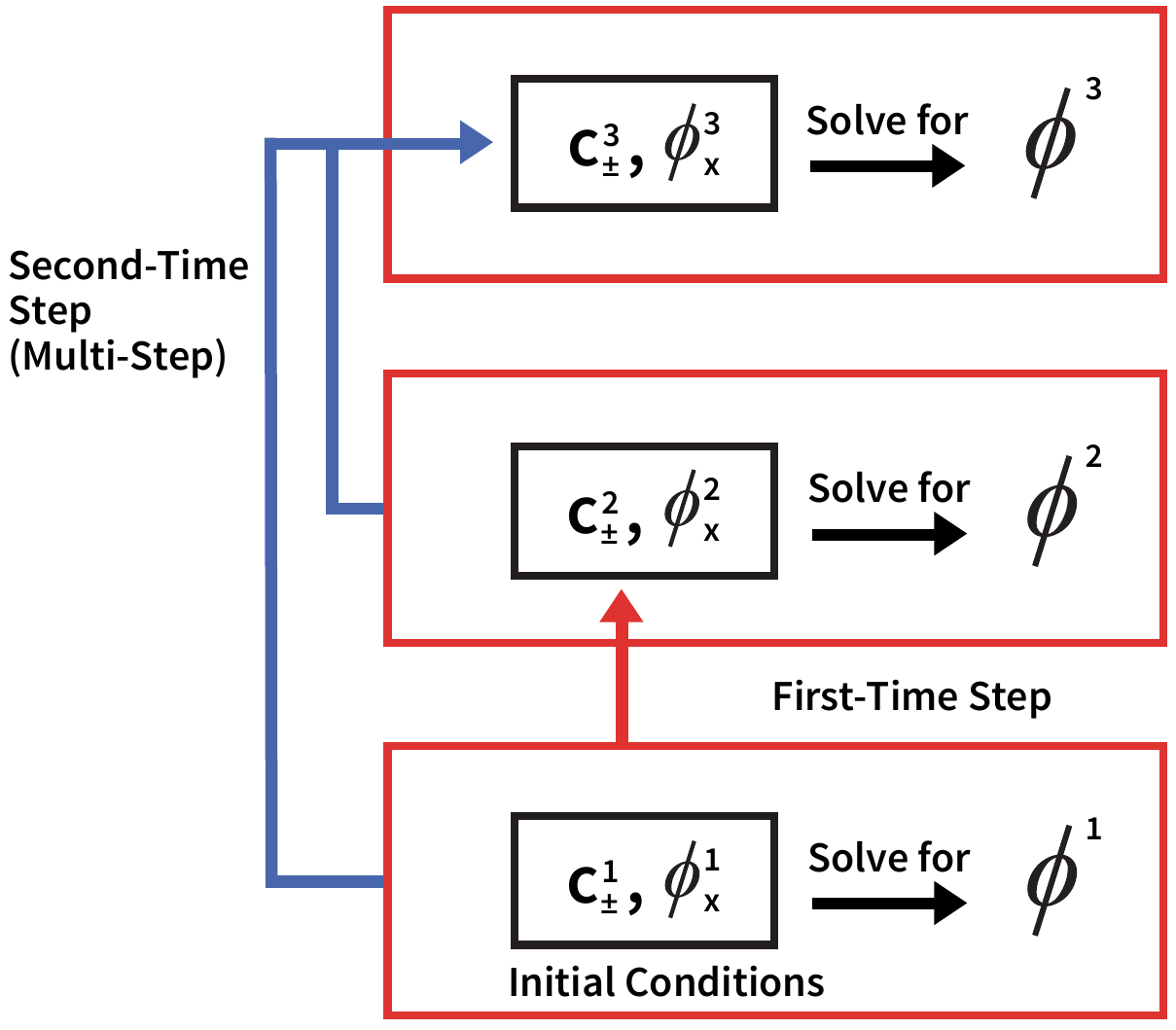}
\caption{The initial
value $\phi_x(1,0)$ is denoted $\phi_x^1$ and $\phi_x^n$ denotes the approximation of $\phi_x(1,t^n)$.  The
initial data $c_\pm(x,0)$ is discretized resulting in initial vectors $\mathbf c_+^1,\, \mathbf c_-^1 \in \mathbb{R}^{N}$.  
Given $\mathbf c^1_\pm$ and $ \phi_x^1$, we solve the elliptic PDE to determine $\pmb{ \phi}^1 \in \mathbb{R}^N$. Using the one-step semi-implicit scheme \eqref{one_step_IMEX}
on the parabolic PDEs (\ref{concentration_nondim}) and the ODE (\ref{current_conservation_nondim}), we determine $\mathbf c_\pm^2$ and $\phi_x^2$. Given $\mathbf c^2_\pm$ and $ \phi_x^2$, we solve the elliptic PDE to determine 
$\pmb{ \phi}^2$.
 \eqref{sbdf2} or \eqref{vssbdf2} is then applied to determine $\mathbf c_\pm^3$ and $\phi_x^3$;
$\pmb{\phi}^3$ is then solved for, using $\mathbf c_\pm^3$ and $ \phi_x^3$.  And so forth.
%
 }
\label{timelevels}
\end{figure}
Figure \ref{timelevels} presents this splitting scheme for the first and second (and subsequent) time-steps. 
It is important 
that $\pmb{\phi}^{n+1}$ be solved for using 
$\mathbf c^{n+1}_\pm$ and $\phi_x^{n+1}$; not doing so reduces the global truncation error from $O(dt^2)$ to $O(dt)$. 
Furthermore, because $\pmb{ \phi}^n$ at time $t^n$ is determined from the other quantities at time $t^n$\textcolor{blue}{,}
the elliptic equation does not factor into the computation of the approximate local truncation error during adaptive 
time-stepping.

\subsubsection{Using the Fully-implicit Schemes BDF2 and VSBDF2}
\label{fully_implicit}

If fully-implicit time-stepping is used, the problem is viewed as a system of DAEs.
If there is an imposed voltage, 
the $3N$ unknowns $\c_+^{n+1}$, $\c_-^{n+1}$, and $\pmb{\phi}^{n+1}$ are 
simultaneously solved for. 
The parabolic PDEs \eqref{concentration_nondim} and 
the Butler-Volmer boundary conditions
\eqref{bv_nondim_L}--\eqref{bv_nondim_R}  are applied at time $t^{n+1}$; this
results in $2N$ equations. 
The elliptic PDE \eqref{poisson_nondim} and its boundary conditions
\eqref{phi_bc_nondim_L}--\eqref{phi_bc_nondim_R} are applied at time
$t^{n+1}$ as well; this provides the additional $N$ equations.

The $3N$ equations are rewritten as $\mathcal{F}(\c_+,\c_-,\pmb{ \phi}) = \vec{0}$
where $\mathcal{F}: \R^{3N}\to\R^{3N}$.  The function $\mathcal{F}$ 
depends on the known quantities $\c_\pm^n$, $\pmb{\phi}^n$, $\c_\pm^{n-1}$, 
and $\pmb{\phi}^{n-1}$.
We use two methods of finding approximate solutions. 
The first is a hand-coded Newton-Raphson method, for which iterations were stopped when $\|\mathcal{F}(\c_+,\c_-,\pmb{ \phi})\| < 10^{-10}$. The second is MATLAB's {\tt fsolve} routine, which, by default, uses a trust region optimization algorithm. The results in Section \ref{method_comparison} show that MATLAB's {\tt fsolve} is the faster and more stable method (as one would expect).

If there is an imposed current, rather than an imposed voltage, 
then 
$\phi_x^{n+1}$ is also solved for, resulting in
 $3N+1$ unknowns.
The additional equation in the 
system of DAEs is \eqref{current_conservation_nondim} .

\section{Simulations of the PNP-FBV Equations}
\label{numerical_tests_PNP}

Table \ref{convergencetest} presents the results from convergence testing performed on the PNP-FBV system using both the SBDF2 and BDF2 numerical schemes.  For these convergence tests, we considered current, rather than voltage boundary conditions, which requires time-stepping the additional ODE \eqref{current_conservation_nondim}. The results demonstrate that both schemes are second-order accurate when the ``ghost point'' implementation of the boundary
conditions \eqref{ghostpoints_left} is used and that using the ``direct'' method
\eqref{bc2_direct} causes a loss of an order of accuracy.
\begin{table}[htb!] 
\centering
\vspace{-.1in}
\begin{tabular}{|c|c|c|c|}
\hline
 $dt$   & SBDF2, direct & SBDF2, ghost & BDF2\\ \hline
 $5 \times 10^{-7}$& 1.9509 & 3.9904 &3.9967 \\ \hline
$5/2\times 10^{-7}$& 1.9758& 3.9952& 3.9983\\ \hline
 $5/4\times 10^{-7}$& 1.9880 & 3.9976&3.9992\\ \hline
$5/8\times 10^{-7}$& 1.9940 &  3.9988& 3.9996 \\ \hline
$5/16 \times 10^{-7}$& 1.9970 & 3.9993   & 3.9998\\ 
\hline
\end{tabular}
\caption{Convergence tests on 
PNP-FBV system \eqref{concentration_nondim}--\eqref{phi_bc_nondim_R} 
with current boundary conditions \eqref{current_conservation_nondim} and an imposed external current
$j_\text{ext}=0.5$. All model parameters were set to $1$ except $\epsilon=0.01$.  
For each method, we compute seven solutions on a uniform mesh $dx = 1/30$ up to time $T = 10^{-5}$, using equation \eqref{one_step_IMEX} for the first time-step of the SBDF2 scheme and equation \eqref{backward_euler} for the first 
time-step of the BDF2 scheme. 
The SBDF2 scheme is implemented using both direct boundary conditions
\eqref{bc2_direct}
and the ghost point approach \eqref{ghostpoints_left}.
The i$^\text{th}$ solution, denoted $\mathbf u_i = [\c_{+,i};\c_{-,i};\phi_{x,i}^{n+1}]$, is computed using $dt=1/2^i \times 10^{-6}$.  
That is $\u_i \in \R^{2 (N+1)+1} \times \R^{n_i}$ where $n_i =  10 \times 2^i+1$ and
$\u_i^{n_i} \in \R^{2 (N+1)+1}$ approximates the solution at the final time $T$.
The discrete solutions $\mathbf u_i$ at the final time $T=10^{-5}$ are used to define the ratios:
$\| \mathbf u^{n_i}_i- \mathbf u^{n_{i+1}}_{i+1} \|/\| \mathbf u^{n_{i+1}}_{i+1}- \mathbf u^{n_{i+2}}_{i+2}\|$ with the $l^2$ norm.  Ratios approaching 4 indicate global second-order accuracy in time. While we only show 5 values of $dt$ in the table, runs with $dt=5/32\times 10^{-7}$ and $dt=5/64\times 10^{-7}$ were necessary to create the ratios.
Richardson extrapolation was not used.
}
\label{convergencetest}
\end{table}

We consider an initial value problem
for the PNP-FBV system
\eqref{concentration_nondim}--\eqref{phi_bc_nondim_R} with constant
imposed voltage.  
Solutions 
are computed using both the (fully-implicit) VSBDF2
adaptive time-stepper and the (semi-implicit) VSSBDF2 adaptive
time-stepper.

The plots to the left in Figure \ref{Full_Model_tmax_effect} present
results from the simulations using the
 (fully-implicit) VSBDF2 adaptive time-stepper.  The bottom plot demonstrates that the
time-step increases 
until it reaches the user-specified
$dt_\text{max}$; the time-step remains at this value for the duration
of the simulation.  The middle plot demonstrates that the
(approximate) local truncation error stays within the user-specified
range of $(tol-range,tol+range)$ until the time when the time-step reaches
$dt_\text{max}$.  Once the time-step stays at $dt_\text{max}$, the
error control mechanism is no longer in play --- the local truncation 
error begins to decay to around $10^{-11}$ and then stays around
this value.  
Similarly, the top plot demonstrates that the norms of the deviation
from the steady-state solution, $| \c^n_\pm - \c_{ss,\pm}|$ and
$|\pmb{ \phi}^n - \pmb{\phi}_{ss}|$, initially decreases exponentially fast, lingers around
$10^{-7}$ for a bit, and then decreases to around $10^{-11}$.  The deviations
don't decrease to $10^{-14}$ as they did for the SBDF2 simulation because
the VSBDF2 adaptive time-stepper uses an iterative nonlinear solver and
this solver is only solving the equations up to a tolerance
of about $10^{-11}$.

The behaviour observed in the plots on the left is what one would expect
for a numerically stable scheme that is computing a solution that is converging
to an asymptotically stable steady-state solution.  This expected
behaviour is {\it not} what is observed
when using the (semi-implicit) VSSBDF2 adaptive time-stepper.

The top figure on the right demonstrates that
the solution 
found by the  (semi-implicit) VSSBDF2 adaptive time-stepper 
initially decays 
towards the numerical steady-state solution.  However, once the
solution is within (approximately) $10^{-7}$ of the steady-state
solution, this convergence ends and the computed solution stays about
$10^{-7}$ away from the steady-state solution.   (We note that if we decrease $tol$ 
then the deviation decreases
further before ``stalling out''.) The middle figure on the right
demonstrates
that the VSSBDF2 adaptive time-stepper is keeping the (approximate) local
truncation error \eqref{final_error}
within the user-specified range of $(tol-range,tol+range)$, as required as long
as $dt_\text{min}< dt < dt_\text{max}$.
The bottom figure on the right
demonstrates that the time-step size initially increases exponentially
fast and after a while
it stabilizes; we denote the value that it
stabilizes to as $dt_\infty$.  To approximate it, we average 100 sequential time-steps at the end of the period of constant voltage/current.

The dashed line in the bottom figure 
on the left is
the stability restriction found by the linear stability analysis
presented in Section 5
of the companion article \cite{YPD_Part2}; we denote this value as
$dt^*$.  This simulation demonstrates that the VSSBDF2 adaptive time-stepper
appears to stabilize at a time-step size $dt_\infty$ that is precisely
the stability restriction $dt^*$.

An example of an adaptive SBDF2 time-stepping scheme can be found in 
Rosam, Jimack and Mullis \cite{Rosam2007} use an adaptive SBDF2 time-stepping scheme to study a problem in binary alloy solidification. In their Figure 4, they appear to show time-steps
relaxing to a constant value (i.e. stabilizing), but its cause is
not given: they report that it is related to the tolerance set in the adaptive 
time-stepper. We find that our $dt_\infty$ is not affected by the tolerance
parameter $tol$.



When using the (fully-implicit) VSBDF2 scheme, 
one can either extrapolate boundary conditions forward in time or 
they can be part of the 
nonlinear DAE system that has to be solved.  We
find that if we extrapolate the coupled boundary conditions of the elliptic PDE then this introduces a stability restriction on the time step size.  
That is, 
the observed behaviour of the solution is no longer like that shown 
in the left plots of Figure \ref{Full_Model_tmax_effect}, rather it becomes like that
shown in the right plots.  
This is in contrast to \cite{Bruno_Jimenez_2014}, where the authors time step a nonlinear diffusion equation using BDF methods. However they extrapolate the diffusivity coefficient forward in time to avoid the need for a nonlinear solve, but find that their time-stepper remains unconditionally stable. 

\subsection{Parameter choice and observed numerical behaviour}

In all of our simulations using the (semi-implicit) VSSBDF2 adaptive time-stepper,
we find that if there was a sufficiently long period of time during which either the 
voltage or the current were held constant then the solution would ``come close to 
but not converge to'' the corresponding steady-state solution (as in the top right plot of 
Figure \ref{Full_Model_tmax_effect})
and the time-step sizes would 
stabilize to a limiting value $dt_\infty$ (as in the bottom right plot of Figure \ref{Full_Model_tmax_effect}).  We find that
$dt_\infty$ depends on the imposed constant voltage (or on the imposed constant current) and that it does not depend on
the parameters $tol$ and $range$ in the adaptive time-stepper.

In the companion article \cite{YPD_Part2} we perform a linear stability 
analysis of the SBDF2 scheme about the steady-state solution, and demonstrate
that the scheme is conditionally stable with a step size stability restriction $dt^*$. 
As suggested in the bottom right plot of Figure \ref{Full_Model_tmax_effect}, and more thoroughly studied
in  \cite{YPD_Part2}, the $dt_\infty$ found by the VSSBDF2 adaptive time-stepper
coincides with $dt^*$.

When computing solutions of the PNP-FBV system with time-independent voltage or current, we find that if we use
the SBDF2 time-stepper with (fixed) $dt > dt^*$ then, after a transient, the computed solutions grow exponentially fast.
If we use a (fixed) $dt < dt^*$ then, after a transient, 
the solution eventually converges exponentially fast to the steady-state solution and the local truncation error decays exponentially fast to zero (round-off).   This is how
we compute the steady-state solution that is used in computing the 
deviations from the steady-state solution presented in Figure \ref{Full_Model_tmax_effect}.

If we use the (semi-implicit) VSSBDF2 adaptive time-stepper then one of two things happens.  If $dt_\text{max} < dt^*$ then as the solution starts to converge to the steady-state solution $dt$ increases to $dt_\text{max}$.  After this happens the local truncation error converges to zero and the solution converges to the steady-state solution.  Alternatively, if $dt_\text{max} \geq dt^*$ then the solution
starts to converge to the steady-state solution.  As a result of this, the
time-steps need grow in order to keep the local truncation error greater than 
$tol-range$ (which has been chosen to be positive). 
Once the time-step size increases past $dt^*$, the unstable eigenmode(s) of the 
underlying linearized problem
start to grow exponentially and this growth stops the solution from converging to the 
steady-state solution.   This growth causes the local truncation error to exceed
$tol+range$, causing the adaptive time-stepper to decrease the time-step size.  The
net effect is that the adaptive time-stepper stabilizes its time-steps at $dt^*$ in order
to keep the local truncation error in the user-specified interval  \cite{YPD_Part2}.

\subsection{Response to Time-dependent Forcing}
\begin{figure}[htb!]
\centering
\includegraphics[width=1\linewidth,height=0.7\linewidth]{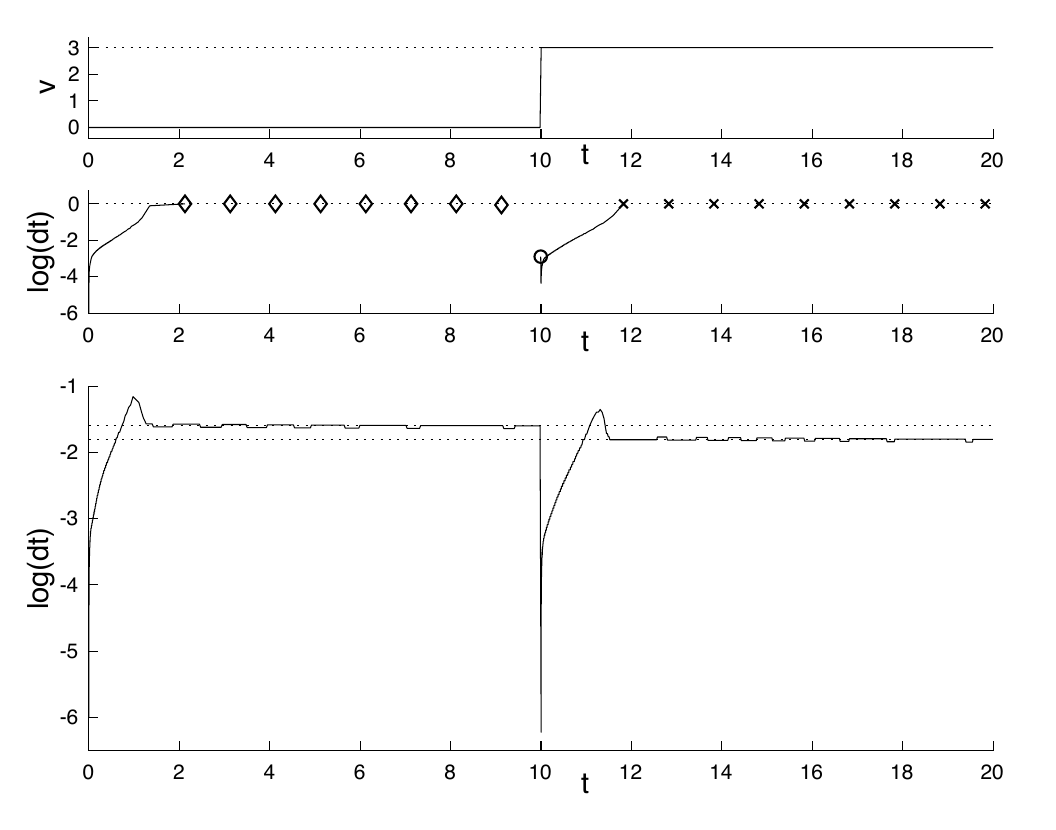}
\caption{PNP-FBV system \eqref{concentration_nondim}--\eqref{phi_bc_nondim_R} with all physical parameters set to $1$ except $\epsilon=0.5$.
A uniform mesh is used: $N = 90$ and $L_1 = L_2 = L_3 = 1/3$.
The initial data for the concentrations is $c_\pm(x,0)=1 + .1 \, \sin(2 \pi x)$. \underline{Top plot:} The imposed time-dependent voltage $v(t)$ has a sharp transition from 0 to 3 at time $t=10$, modeled by the equation $v(t) = 3(\tanh(1000(t-10))+1)/2$. 
\underline{Middle plot:} The time steps as a function of time for the fully-implicit scheme (solid line) and the semi-implicit scheme (dashed line).  The diamond and x-marks indicate time-steps that are within 95\% of the user-specified maximum value $dt_\text{max}$.  The circle denotes the first time-step size after the last diamond-marked time step.
\underline{Bottom plot:} The time steps as a function of time for the semi-implicit scheme (solid line) and the fully-implicit scheme (dashed line). The dotted lines in the center and right figures correspond to $dt^* = 0.0250$ and $0.0158$; they are the stability restrictions found using the linear stability analysis presented in Section 4 of the companion paper. }
\label{dt_vs_time}
%
\end{figure}

We next considered the PNP-FBV system with time-dependent imposed voltage.  The top plot of Figure \ref{dt_vs_time} shows
the imposed voltage: essentially a step function with a fast transition at $t=10$.

%

The middle plot of Figure \ref{dt_vs_time} presents the time-steps chosen by the
(fully-implicit) VSBDF2 adaptive time-stepper.  In the first half of the simulation, the
solution is converging to the $v=0$ steady-state solution.  
During this time, the time-steps coarsen until they reach 
the user-specified maximum time-step size: $dt_\text{max}=1$.  See the diamond-marked times. The time-step refines 
markedly in response to the change in imposed voltage; the open circle indicates the 
first time-step size after the last diamond-marked time-step size.  After the voltage 
step, the
solution is converging to the $v=3$ steady-state solution.  Again, the time-steps coarsen until they reach $dt_\text{max}$ (see the x-marks in the plot).
The deviation of the solution from the steady-state solution and the behaviour 
of the local truncation error are the analogues of that shown in plots on the left of
Figure \ref{Full_Model_tmax_effect}. 

The bottom plot of Figure \ref{dt_vs_time} presents the time-steps chosen by the
(semi-implicit) VSSBDF2 adaptive time-stepper.   As discussed for the plots in the
left of Figure \ref{Full_Model_tmax_effect}, the solution starts to converge to the
$v=0$ steady-state solution and the time-step size start to coarsen,
stabilizing to a value $dt_\infty$.  The time-step refine in response to the change
in the imposed voltage after which the solution starts to converge to the
$v=3$ steady-state solution and the time-step size start to coarsen,
stabilizing to a different value $dt_\infty$.   The two values of $dt_\infty$ coincide
with the stability restrictions $dt^*$ found for the $v=0$ and $v=3$
steady-state solutions \cite{YPD_Part2}.

\subsection{Effect of the Singular Perturbation Parameter}
\begin{figure}[htb!]
\centering
    \begin{subfigure}[htb!]{0.49\textwidth}
        \centering
\caption{}
\includegraphics[width=\linewidth]{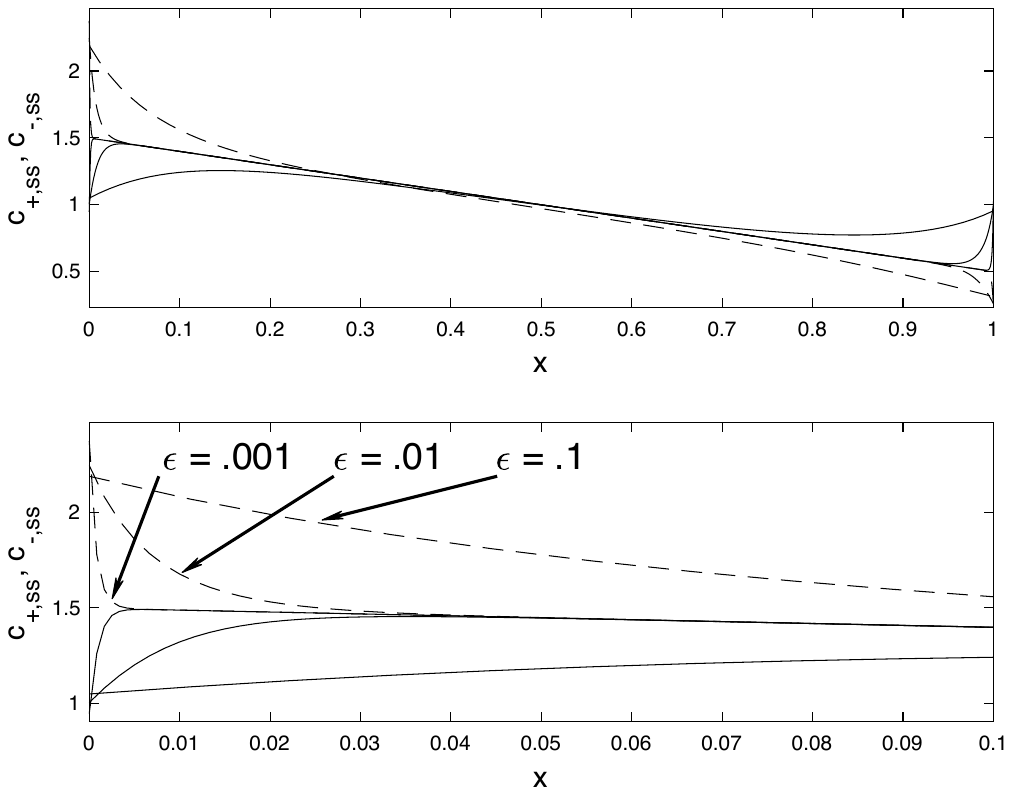}
    \end{subfigure}
    \begin{subfigure}[htb!]{0.49\textwidth}
        \centering
\caption{}
\includegraphics[width=\linewidth]{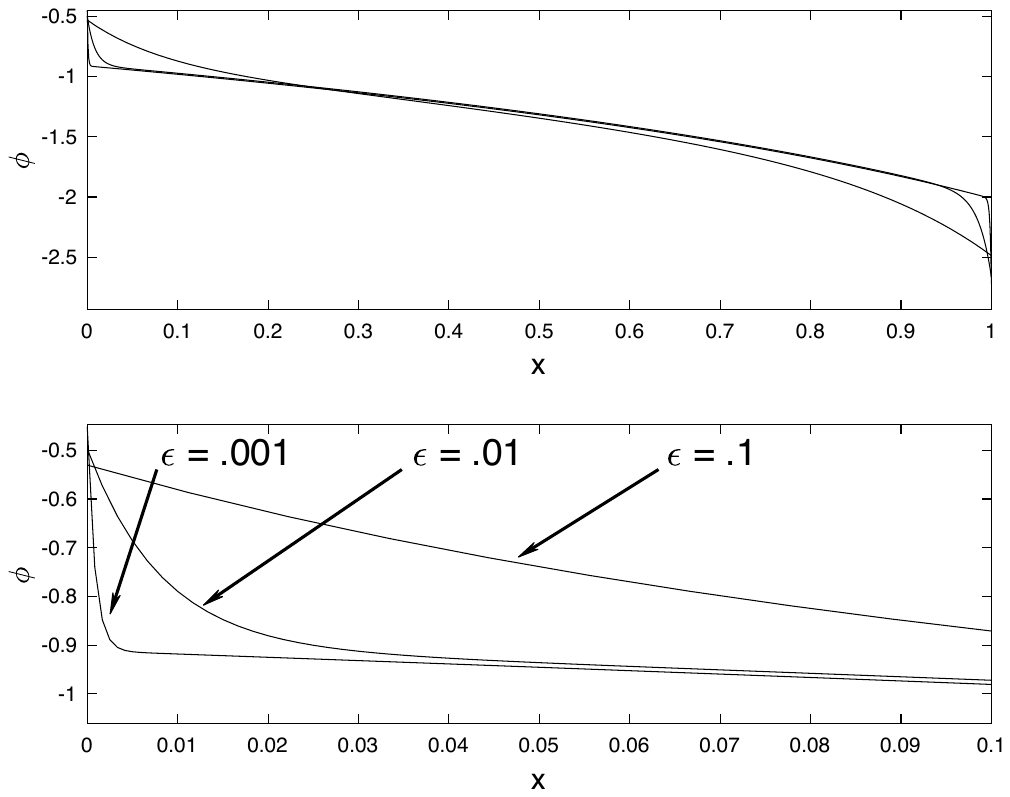}
    \end{subfigure}
\caption{Steady-state solutions for
PNP-FBV system \eqref{concentration_nondim}--\eqref{phi_bc_nondim_R}
with current boundary conditions \eqref{current_conservation_nondim}.
The imposed external current is constant $j_\text{ext}(t)=0.5$.  The simulations 
were run until time $t=5$; the profiles at that time are shown here.  
All physical parameters were set to $1$, except for $\epsilon$.
For $\epsilon = 0.1$, a uniform-in-space mesh was used with 
$dx = 1/90$.  
For $\epsilon = 0.01$, a piecewise uniform-in-space
mesh was used with $dx = 1/600$ in $[0,0.1]$ and $[0.9,1]$ and
$dx = 1/75$ in $[0.1,0.9]$. 
For $\epsilon = 0.001$, a piecewise uniform-in-space
mesh was used with $dx = 1/1200$ in $[0,0.05]$ and $[0.95,1]$ and
$dx = 3/100$ in $[0.05,0.95]$.  \underline{Left plot:} concentrations  $c_{+,ss}$ (solid) and $c_{-,ss}$ (dashed), versus $x$. 
\underline{Right plot:} potential, $\phi_{ss}$, versus $x$.}
\label{SS_solns_vs_eps}
\end{figure}

Next, we consider the singular perturbation
parameter, $\epsilon$, which controls the dynamics of the boundary layers
near the electrodes in the PNP-FBV system \eqref{poisson_nondim}.  
Figure \ref{SS_solns_vs_eps} presents some steady-state solutions for
PNP-FBV system \eqref{concentration_nondim}--\eqref{phi_bc_nondim_R}
with current boundary conditions \eqref{current_conservation_nondim} as computed using
the SBDF2 scheme.
The plots on the left present the steady-state concentrations, $c_{+,ss}(x)$ and $c_{-,ss}(x)$, for three values of $\epsilon$.  The 
smaller the value of $\epsilon$, the thinner the transition layer
near $x=0$ and $x=1$ is.  The plots on the right present the corresponding steady-state 
potential, $\phi_{ss}(x)$.  The top left plot shows that, in the bulk,
$c_{+,ss}$ and $c_{-,ss}$ are very close in value (the system is nearly electroneutral).  We find that, for small
values of $\epsilon$,
$c_{+,ss}(x) - c_{-,ss}(x) \approx \epsilon^2 \, C(x)$ in the bulk.  This is the reason why the
steady-state potential is not linear in the bulk (see the top right plot); in the bulk, it satisfies $\phi_{xx} = C/2$.

%
%
%
%

For relatively small values of $\epsilon$, we find that if one decreases $\epsilon$ by a factor of $a$ then the limiting
time-step size of the (semi-implicit) VSSBDF2 adaptive time-stepper, $dt_\infty$, decreases by roughly $a^2$; see Table \ref{epsilon_test}. A detailed discussion of this effect can be found in the companion article \cite{YPD_Part2}.
\begin{table}[htb!] 
\centering
\begin{tabular}{|l|l|l|l|}
\hline
                     $\epsilon$   & $dt_\infty$ & $dt_\infty/\epsilon^2$  \\ \hline
                      $0.1$ &  $0.0257$ & $2.57$ \\ \hline
                       $0.01$& $0.000241$ & $2.41$ \\ \hline
                      $0.001$&  $3.70 \times 10^{-6}$& $3.70$ \\ \hline
                       $10^{-4}$&  $2.48 \times 10^{-8}$ & $2.48$ \\ \hline
                       $10^{-5}$&  $2.25 \times 10^{-10}$ & $2.25$\\ \hline
                       $10^{-6}$& $2.37\times 10^{-12}$ & $2.37$ \\
\hline
\end{tabular}
\caption{The PNP-FBV system \eqref{concentration_nondim}--\eqref{phi_bc_nondim_R} are simulated with the
imposed voltage $v(t)=0$.   All physical parameters were set to $1$, except for 
$\epsilon$.
Six values of $\epsilon$ are considered. We denote the limiting timestep size found by the VSSBDF2 adaptive time-stepper by $dt_\infty$. $dt_\infty$ was computed by
averaging the time-step sizes over 100 sequential time-steps chosen late in the simulation.}
\label{epsilon_test}
\end{table}

By taking small enough time-steps, the (semi-implicit) VSSBDF2 adaptive 
time-stepper is
able to run with small 
values of $\epsilon$. Furthermore, although it is more computationally expensive, the (fully-implicit) VSBDF2 adaptive time-stepper is able to take arbitrarily large time-steps for equilibrated solutions, regardless of the value of $\epsilon$.

\subsection{Comparison Between Semi-implicit and Fully-Implicit Methods}
\label{method_comparison}
As demonstrated in Table \ref{epsilon_test}, the stability restriction on the time-step size for VSSBDF2 is roughly on the order of $\epsilon^2$ as $\epsilon \to 0$, whereas there is no dependence on $\epsilon$ on the step size for VSBDF2. This implies that there will be a break-even point: in some regimes the VSSBDF2 scheme will be faster, while in other regimes the VSBDF2 scheme will be faster despite being significantly slower per time-step. 

To compare the two schemes, we repeated the experiments in Table \ref{epsilon_test},
using Matlab's profiler to measure performance time, and keeping track of the number of iterations and the average time step size for our simulations.  See Table \ref{profile_table}. ``Total iterations" counts the number of coarse-step-fine-step sequences taken, including those taken while seeking an acceptable time-step size (see Algorithm \ref{adaptive_algorithm}). ``Time steps" counts the number of times the time index is actually incremented. When using VSBDF2, we tested two methods of solving the system of nonlinear equations: Newton's method, and Matlab's \texttt{fsolve}. We ran all simulations to a final time $t_{\text{end}}=1$. We did not run with values of $\epsilon$ smaller than $10^{-4}$ because the VSSBDF2 code would have taken over fifteen hours.  Also, our Newton-Raphson method code ran into conditioning issues for $\epsilon = 10^{-4}$ although Matlab's \texttt{fsolve}
did not.

\begin{table}[]
\begin{tabular}{|l|l|l|l|l|}
\hline

$\epsilon$ & Performance time (s)  & \begin{tabular}[c]{@{}l@{}}Total iterations (including \\ error control iterations)\end{tabular} & Time steps & Average time step \\ \hline
        & \multicolumn{4}{c|}{\textbf{(semi-implicit) VSSBDF2}}                                                                                                                           \\ \hline
0.1     & 0.661                 & 625                                                                                                 & 445        & 2.21E-03          \\ \hline
0.01    & 5.524                 & 7836                                                                                                & 7614       & 1.31E-04         \\ \hline
0.001   & 552.675               & 788026                                                                                              & 754002     & 1.23E-06          \\ \hline
10$^{-4}$    & \textgreater 15 hours &                                                                                                     &            &                   \\ \hline
        & \multicolumn{4}{c|}{\textbf{(fully-implicit) VSBDF2, Newton solver}}                                                                                                         \\ \hline
0.1     & 201.201               & 530                                                                                                 & 356        & 3.15E-03          \\ \hline
0.01    & 214.671               & 525                                                                                                 & 354        & 2.86E-03          \\ \hline
0.001   & 255.329               & 530                                                                                                 & 355        & 3.16E-03         \\ \hline
10$^{-4}$    & \multicolumn{4}{l|}{Unable to run, Jacobian poorly conditioned}                                                                                              \\ \hline
        & \multicolumn{4}{c|}{\textbf{(fully-implicit) VSBDF2, \texttt{fsolve} solver}}                                                                                                \\ \hline
0.1     & 157.09                & 530                                                                                                 & 356        & 3.17E-03          \\ \hline
0.01    & 153.099               & 525                                                                                                 & 354        & 2.86E-03          \\ \hline
0.001   & 158.584               & 530                                                                                                 & 355        & 3.16E-03          \\ \hline
10$^{-4}$    & 157.031               & 526                                                                                                 & 355        & 2.86E-03          \\ \hline
\end{tabular}
\caption{Profiling results on the (semi-implicit) VSSBDF2 and (fully-implicit)
VSBDF2 adaptive time-steppers. For the fully-implicit method two nonlinear solvers were compared: a simple Newton method and Matlab's \texttt{fsolve} routine. Parameter values are identical to the ones used for Table \ref{epsilon_test}, except all experiments were run to a final time $t_\text{end}=1$. Performance time was obtained using Matlab's profiler. ``Total iterations" counts number of coarse-step-fine-step sequences taken, including those taken while seeking
an acceptable time-step size (see Algorithm \ref{adaptive_algorithm}), whereas ``Time steps" counts the number of times the time index is actually incremented. }
\label{profile_table}
\end{table}

The profiling results show a clear trend: with decreasing $\epsilon$, the semi-implicit scheme's number of time steps increases rapidly, while the fully-implicit scheme's number of time steps remains constant. This is due to the semi-implicit scheme's time step stability restriction, which depends on $\epsilon$. This shows that, for this particular set of parameters and with no forcing, the semi-implicit scheme is favorable when $\epsilon$ is roughly greater than $0.001$, and the fully-implicit scheme is favorable when $\epsilon < 0.001$. In other parameter and forcing regimes, the ``break-even'' value of $\epsilon$ will be different.
We note that the test has been run in a situation where the solution
equilibrates.  One does this when seeking steady-state solutions, for example.  
The final time $t_\text{end}$ is another user-specified parameter that affects
the ``break-even'' value of $\epsilon$.

\section{Conclusions and Future Work}

In this work, we considered the Poisson-Nernst-Planck equations with generalized Frumkin-Butler-Volmer reaction kinetics at the electrodes. We developed a solver that dynamically chooses the time-step size so that the local truncation error is within user-specified bounds; this is a novel contribution to the solution of the PNP equations. The spatial discretization can be nonuniform, allowing for the finer mesh needed
near the electrodes to resolve the boundary layers. Adaptive time-stepping allows
the solver to choose small time-steps during the initial transient period during which the boundary layers may be forming quickly in response to the boundary conditions, and also allows one to more efficiently study physical situations in which the imposed voltage or current have sudden, fast changes.

We considered two adaptive time-stepping schemes in this work. The (semi-implicit)
VSSBDF2 adaptive time-stepper was shown to have a stability restriction which causes its time-step sizes to stabilize to a limiting value, $dt_\infty$, in the long-time limit of constant voltage or current, whereas the (fully-implicit) VSBDF2 adaptive time-stepper was found to have no restriction on the time step beyond the user-specified $dt_\text{max}$. However, the fully-implicit time-stepper requires a computationally expensive nonlinear solve per step, resulting in longer computation time per time step. 
The cause of the limiting time-step, $dt_\infty$, found by the 
(semi-implicit) VSSBDF2 adaptive time-stepper 
can be understood by linearizing the SBDF2 scheme about the steady-state solution 
and is fully explored in the companion article \cite{YPD_Part2}.

We found that $dt_\infty$ is $\epsilon$-dependent.  As a result, the
(semi-implicit) VSSBDF2 adaptive time-stepper is favoured for larger values
of $\epsilon$ and the (fully-implicit) VSBDF2 adaptive time-stepper is favoured
for smaller values.

\appendix
\section{Derivation of the Local Truncation Error Formula and the Extrapolation Formula}
\label{app:extrapolation}

In this appendix, the error approximation and extrapolation formula used in Subsection \ref{time_stepping_error_control} are derived. We consider time-stepping the ODE
$
{du}/{dt}=f(u)+g(u)$
where the value of $u$ at time $n+1$ is given by \eqref{vssbdf2}.
Making the substitutions $u^{n+1}=u(t+dt_\text{now})$ and $u^{n-1}=u(t-dt_\text{old})$ and performing a multivariable Taylor expansion about $t$, the local truncation error is found to take the form  
\begin{align}
\label{LTE_1}
LTE&=\left(u'(t)-f(u(t))-g(u(t))\right)dt_\text{now} + \left(u''(t)-f'(u(t)) \, u'(t)-g'(u(t)) \, u'(t)  \right) dt_\text{now}^2 \nonumber \\
&\hspace{.5in}+ \left(\frac{1}{2} f''(u(t)) \, u'(t)^2  +\frac{1}{2}f'(u(t)) \, u''(t)  - \frac{1}{6}u^{(3)}(t)) \right) dt_\text{old} \, dt_\text{now}^2 \nonumber \\
&\hspace{0.5in}+ \left(-\frac{1}{2} g''(u(t)) \, u'(t)^2  -\frac{1}{2}g'(u(t)) \, u''(t)+\frac{1}{3}u^{(3)}(t)\right)dt_\text{now}^3 + O(dt^4)
\end{align}
where $O(dt^4)$ denotes terms where $dt_\text{now}$ and $dt_\text{old}$ combined appear four or more times. If $u(t)$ 
is a thrice-differentiable solution of the ODE, the $dt_\text{now}$ and $dt_\text{now}^2$ terms vanish
and the cubic term can be simplified:
\begin{align}
\label{LTE_2}
LTE&=\left(\frac{1}{2}f''(u(t)) \, u'(t)^2  +\frac{1}{2}f'(u(t)) \, u''(t) - \frac{1}{6}u^{(3)}(t)\right)dt_\text{old}dt_\text{now}^2  \nonumber \\
&\hspace{.75in}+ \left(-\frac{1}{2}g''(u(t)) \, u'(t)^2  -\frac{1}{2}g'(u(t))u''(t)+\frac{1}{3}u^{(3)}(t)\right)dt_\text{now}^3 + O(dt^4) \\
\label{LTE_3_app}
&=\left(\frac{1}{3} u^{(3)}(t) - \frac{1}{2} g''(u(t)) \, u'(t)^2- \frac{1}{2} g'(u(t)) \, u''(t) \right)dt_\text{now}^2\left(dt_\text{now}+dt_\text{old}\right) + O(dt^4).
\end{align}
This means that if $u^{(3)}(t)$ is bounded, 
the local truncation error for the coarse step is
\begin{equation}
\label{coarse_error_app}
\epsilon_c=u_c^{n+1}-u(t^{n+1})\approx Cdt_\text{now}^2\left(dt_\text{now}+dt_\text{old}\right) 
\end{equation}
To find the local truncation error for $u_f^{n+1}$, we find the local truncation error for $u_f^{n+\frac{1}{2}}$ by replacing $dt_\text{now}$ and $dt_\text{old}$ in \eqref{coarse_error_app} by $\frac{dt_\text{now}}{2}$ and $\frac{dt_\text{old}}{2}$, respectively. We then find the local truncation error for $u_f^{n+1}$ by setting both $dt_\text{now}$ and $dt_\text{old}$ to $\frac{dt_\text{now}}{2}$. Adding these local truncation errors yields the LTE for $u_f^{n+1}$,
\begin{equation}
\label{fine_error_app}
\epsilon_f=u_f^{n+1}-u(t^{n+1})\approx \frac{Cdt_\text{now}^2\left(dt_\text{now}+dt_\text{old}\right)}{8} + C \, \frac{dt_\text{now}^3}{4} 
\end{equation}
Using equations \eqref{coarse_error_app} and \eqref{fine_error_app}, we can approximate $C$ using $u_c$, $u_f$, $dt_\text{now}$ and $dt_\text{old}$ and then use this approximation in
\eqref{coarse_error_app}, resulting in
\begin{equation}
\label{error}
\epsilon_c\approx\frac{u_c-u_f}{\frac{5}{8}dt_\text{now}^3+\frac{7}{8}dt_\text{now}^2dt_\text{old}} \; dt_\text{now}^2 \left(dt_\text{now}+dt_\text{old}\right)
\end{equation}
Finally, we subtract $\epsilon_c$ from $u_c^{n+1}$ to create a more accurate approximation of $u(t^{n+1})$. This results in
the extrapolation formula
$u(t^{n+1})\approx u_c^{n+1}-\epsilon_c=\alpha u^{n+1}_c + \beta u^{n+1}_f$
where
\begin{equation}
\label{alphabeta_app}
\alpha=-\frac{dt_\text{old}+3 \, dt_\text{now}}{7 \, dt_\text{old}+5 \, dt_\text{now}},\;\;\;\beta=8\, \frac{dt_\text{old}+dt_\text{now}}{7 \, dt_\text{old}+5 \, dt_\text{now}}.
\end{equation}

We note that if the time-stepping scheme is fully-implicit, this corresponds to setting $f(u)=0$ in equation \eqref{alphabeta_app}. In this case, \eqref{LTE_3_app} is replaced by
$$
LTE = - \frac{1}{6} \left( g''(u(t)) \, u'(t)^2+  g'(u(t)) \, u''(t) \right)dt_\text{now}^2\left(dt_\text{now}+dt_\text{old}\right) + O(dt^4)
$$
but there are no resulting changes to equations \eqref{coarse_error_app}--\eqref{alphabeta_app}.

\section{Acknowledgements} 

Research supported in part by NSERC grant OGP06617.

We thank Greg Lewis, Steven Ruuth, and Adam Stinchcombe for helpful conversations and encouragement.



\end{document}